\documentclass{article}

\usepackage{latexsym}
\usepackage[T1]{fontenc}
\usepackage[utf8]{inputenc}
\usepackage[english]{babel}
\usepackage{algorithmic}
\usepackage{xcolor}
\usepackage{mathrsfs}
\usepackage{yfonts}
\usepackage{bm}
\usepackage{amssymb}
\usepackage{amsmath}
\usepackage{amsthm}
\usepackage{mathdots}
\usepackage{graphicx}
\usepackage{tikz-cd}
\usepackage{tikz}
\usetikzlibrary{arrows}
\usepackage{amsfonts}
\newtheorem{Thm}{Theorem}

\newtheorem{Ex}{Example}
\newtheorem{Pro}{Property}
\newtheorem{Lemma}{Lemma}
\newtheorem{Cor}{Corollary}
\newtheorem{Remark}{Remark}

\usepackage{booktabs}
\usepackage{array}
\newcolumntype{C}{>{\quad}c<{\quad}}

\usepackage{graphicx}
\usepackage{epsfig}
\usepackage{framed}
\usepackage{mathrsfs}
\usepackage{mathtools}

\usetikzlibrary{arrows, positioning, automata}
\usepackage{ytableau}
\ytableausetup
{boxsize=1.1em}
\usepackage [ruled,linesnumbered]{algorithm2e}
\SetAlgoCaptionSeparator{.}
\SetAlgoCaptionLayout{small}
\SetKw{KwGo}{go to}
\SetKw{KwBreak}{break}

\usepackage{subfigure}
\usepackage{multirow}

\usepackage{appendix}

\title{On a conjecture on prime double square tiles}

\author{Michela Ascolese, Andrea Frosini\\ Dipartimento di Matematica e Informatica \\ Università degli Studi di Firenze \\ email: \{michela.ascolese, andrea.frosini\}@unifi.it}

\begin{document}
\maketitle
\begin{abstract}
In \cite{BGL}, while studying a relevant class of polyominoes that tile the plane by translation, i.e., double square polyominoes, the authors found that their boundary words, encoded by the Freeman chain coding on a four letters alphabet, have specific interesting properties that involve notions of combinatorics on words such as palindromicity, periodicity and symmetry. Furthermore, they defined a notion of reducibility on double squares using homologous morphisms, so leading to a set of irreducible tile elements called \emph{prime double squares}. The authors, by inspecting the boundary words of the smallest prime double squares, conjectured the strong property that no runs of two (or more) consecutive equal letters are present there. In this paper, we prove such a conjecture using combinatorics on words tools, and setting the path to the definition of a fast generation algorithm and to the possibility of enumerating the elements of this class w.r.t. standard parameters, as perimeter and area.  

\end{abstract}

\section{Introduction}
Algorithmic studies of planar tilings greatly benefit from the seminal works \cite{Wang} and \cite{WL}, where the decidability of the existence of planar tilings is addressed both with a given set of tiles and with a single one. In the first case, it has been shown that each Turing machine computation can be simulated by a planar tiling using a suitable set of tiles (the author used colored dominoes) without rotations, from an initial partial configuration that models the input tape of the Turing machine. On the other hand, if one only tile is provided, in particular if an {\em exact polyomino} tile is considered, i.e., a four connected finite set of points in $\mathbb{Z}^2$ that tiles the plane by translation, the computation becomes much easier. As a matter of fact in \cite{WL} the authors, with the aim of proving a conjecture by Shapiro, polynomially detected exact polyominoes by surrounding each of them with four or six copies of itself. Moving from this result, Beauquier and Nivat in \cite{BN} characterized the border of an exact polyomino, regarded as a word on a four letters alphabet through the Freeman chain coding, using the notions of rotation and conjugation, proper of combinatorics on words, and setting a strong connection between these two areas. In particular, the boundary word $P$ of such an exact polyomino can be factorized according to the equation $P=X_{1}X_{2}X_{3}\widehat{X}_{1}\widehat{X}_{2}\widehat{X}_{3}$, where $\widehat{X}$ refers to the word $X$ considered as a path and travelled in the opposite direction. According to \cite{BN}, at most one among $X_{1}$, $X_{2}$ and $X_{3}$ can be empty, and we refer to pseudo-squares in this case, pseudo-hexagons otherwise.  

\begin{figure}[ht]
        \centering
        \includegraphics[height=2.3cm]{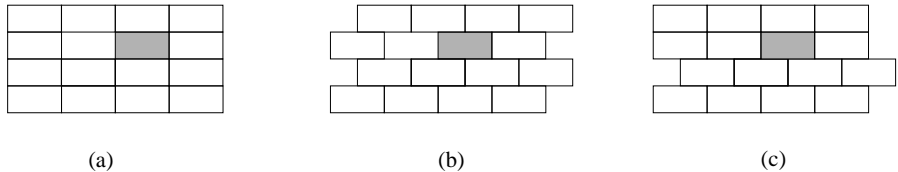}
        \caption{Three tilings of the plane with the two cell domino showing different behaviours. In (a) the domino acts as a pseudo-square, since it is surrounded by four copies of itself. In (b) the domino acts as a pseudo-hexagon, since six copies surround each domino. In (c) both the behaviours are present. Here, the dark domino is surrounded by five copies of itself, but this is not the case for each element of the tiling.}
        \label{fig:tilings}
    \end{figure}

It is easy to verify that an exact polyomino $P$ can tile the plane in different ways, and that it can show both pseudo-hexagon and pseudo-square behaviour also in the same tiling. However, the arrangement of the copies of $P$ in each of its tilings has a periodical behaviour along one or two discrete directions, that are strictly related to the choice of the decomposition of $P$. See Fig.~\ref{fig:tilings} for examples.

Relying on these strong geometrical properties, in \cite{BBL} it was proved that an exact polyomino tiles the plane as a pseudo-square in at most two distinct ways. Furthermore, if two different pseudo-square factorizations of $P$ exist, then no decomposition as a pseudo-hexagon does. The authors refer to these exact polyominoes as {\it double squares}. In Fig.~\ref{fig:ex}, three double square polyominoes are shown together with their factorizations.

It is shown in \cite{BGL} that double squares have specific border properties that are exploited using equations on words. These properties lead to the definition of two operators that allow to exhaustively generate double squares, and to the notion of {\it double square reducibility} through the definition of suitable {\it homologous morphisms}. Among double squares, those that can be reduced without intermediate steps to the simple square constitute the subclass of {\it prime} double squares, and are pointed out for their relevance in the exhaustive generation of double squares. 

Still in \cite{BGL}, the following conjecture was proposed

\medskip

{\bf Conjecture 35.} Let $w$ be the boundary word of a prime double square tile in a four letters alphabet $\Sigma$. Then, for any letter $\alpha \in \Sigma$, $\alpha \alpha$ is not a subword of $w$. 

\medskip

In this article, we prove the conjecture by providing a strong combinatorial structure of the boundary word of any prime double square. In particular, in the next section we recall basic definitions on combinatorics on words and some preliminary results to approach the conjecture. 

In Section~\ref{sec:properties}, we give some properties of the double squares' boundary word, mainly using the results from \cite{BGL}. The proof of the conjecture is in Section~\ref{sec:conjecture}, while Section~\ref{sec:conclusion} is devoted to point out future combinatorial, algorithmic and geometrical research paths that originate from our main result.

\section{Basic notions and previous results}

In this section we fix the notation
and recall some basic results on exact polyominoes to introduce the study of Conjecture 35.

A {\it polyomino} is a subset of the square grid $\mathbb{Z}^{2}$ whose boundary is a continuous, closed and non-intersecting path. We describe polyominoes by coding their boundary through a word defined on the alphabet $\Sigma=\{0,1,\overline{0},\overline{1}\}$, whose elements correspond to the directions  $\{\rightarrow,\uparrow,\leftarrow,\downarrow\}$ of steps made in the grid, respectively. We say that the letters $0$ and $\overline{0}$, resp. $1$ and $\overline{1}$, are \emph{conjugate}, since they represent opposite directions (see Example~\ref{ex:prime_ex}). We indicate by $\Sigma^{*}$ the free monoid on $\Sigma$, i.e., the set of all words defined on the alphabet $\Sigma$, being $\varepsilon$ the empty one. We therefore call \emph{simple square} the polyomino coded by the word $S=01\overline{0}\overline{1}$. Given a word $w\in\Sigma^{*}$, $|w|$ indicates its length, $|w|_{\alpha}$ indicates the number of occurrences of the letter $\alpha$ in $w$, and $w^{n}$ indicates the concatenation of $n$ copies of the word itself.

Boundary words can be considered as \emph{circular} words, since the coding of the boundary of the polyomino can be defined up to the starting point.
Furthermore, the conditions of closure and non-intersection of a polyomino boundary word $P$ can be stated as $|P|_{\alpha}=|P|_{\overline{\alpha}}$ and $|Q|_{\alpha}\neq |Q|_{\overline{\alpha}}$ for each $Q$ proper subword of $P$, say $Q\subset P$, and $\alpha=0,1$.

We define three operators on a word $w=w_{1}w_{2}\dots w_{n}\in\Sigma^{*}$:
\begin{enumerate}
    \item the \emph{conjugate of $w$}, $\overline{w}$, is the word obtained by replacing each letter of $w$ with its conjugate;
    \item the \emph{reversal of $w$}, $\Tilde{w}$, is defined as $\Tilde{w}=w_{n}w_{n-1}\dots w_{1}$. A \emph{palindrome} is a word s.t. $w=\Tilde{w}$;
    \item the \emph{hat of w}, $\widehat{w}$, is the composition of the previous operations, $\widehat{w}=\Tilde{\overline{w}}$.
\end{enumerate}

Finally, an \emph{exact} polyomino is a polyomino that tiles the plane by translation.
Beauquier and Nivat characterized exact polyominoes in relation to their boundary word, providing the following 

\begin{Thm}(\cite{BN})
A polyomino $P$ is exact if and only if there exist $X_{1},X_{2},X_{3}\in\Sigma^{*}$ such that
$$P=X_{1}X_{2}X_{3}\widehat{X}_1\widehat{X}_2\widehat{X}_3,$$ where at most one of the words is empty. This factorization may be not unique.
\end{Thm}

We will refer to this decomposition as a \emph{BN-factorization}. Starting from their BN-factorization(s), exact polyominoes can be further divided in two classes:
\emph{pseudo-hexagons}, if $X_{1},X_{2}$ and $X_{3}$ are all non-empty words, and \emph{pseudo-squares}, if one of the words is empty. We will focus on \emph{double square polyominoes}, i.e., those ones that admit two different BN-factorizations as a square, $P=AB\widehat{A}\widehat{B}=XY\widehat{X}\widehat{Y}$.
Due to the presence of two BN-factorizations, double squares' boundary words can be written in the general form obtained from Corollary $6$ in \cite{BPF},
\begin{equation}\label{eq:genform}
P=w_{1}w_{2}w_{3}w_{4}w_{5}w_{6}w_{7}w_{8},
\end{equation}
where $A=w_{1}w_{2}$, $B=w_{3}w_{4}$, $\widehat{A}=w_{5}w_{6}$, $\widehat{B}=w_{7}w_{8}$ and $X=w_{2}w_{3}$, $Y=w_{4}w_{5}$, $\widehat{X}=w_{6}w_{7}$, $\widehat{Y}=w_{8}w_{1}$, with $w_1, \dots, w_8$ non empty. 

We now introduce the notion of \emph{homologous morphism}. A morphism is a function $\varphi:\Sigma\rightarrow\Sigma^{*}$ s.t. $\varphi(\alpha\beta)=\varphi(\alpha)\varphi(\beta)$ with $\alpha,\beta\in\Sigma$, i.e., it preserves concatenation, and it is said to be \emph{homologous} if $\varphi(\widehat{A})=\widehat{\varphi(A)}$ for all $A\in\Sigma^{*}$, i.e., it preserves the hat operation. From now on we will refer to homologous morphisms only. For each exact polyomino $P=AB\widehat{A}\widehat{B}$, we can define the trivial morphism that maps the simple square in $P$ as $\varphi_{P}(0)=A$, $\varphi_{P}(1)=B$.
In general, the boundary word of an exact polyomino can be obtained starting from the simple square through the composition of two or more morphisms (see Example~\ref{ex:composizione}). A double square is \emph{prime} if its boundary word $P$ is such that, for any homologous morphism $\varphi$, the equality $P = \varphi(U)$ implies that either $U=P$ or $U$ is the boundary word of the simple square. This property can be rephrased saying that a double square is prime if its trivial morphism  can not be obtained by composing two or more different morphisms.

\begin{Ex}\label{ex:prime_ex}
In Fig.~\ref{fig:ex} are shown three prime double squares, whose boundary words are, from left to right, $P_{(a)}=1\overline{0}1\overline{0}101010\overline{1}0\overline{1}0\overline{1}\overline{0}\overline{1}\overline{0}\overline{1}\overline{0}$, $P_{(b)}=1\overline{0}1010\overline{1}0\overline{1}\overline{0}\overline{1}\overline{0}$ and $P_{(c)}=1\overline{0}1\overline{0}1010\overline{1}0\overline{1}0\overline{1}\overline{0}\overline{1}\overline{0}$.

\begin{figure}[ht]
    \centering
    \includegraphics[height=4cm]{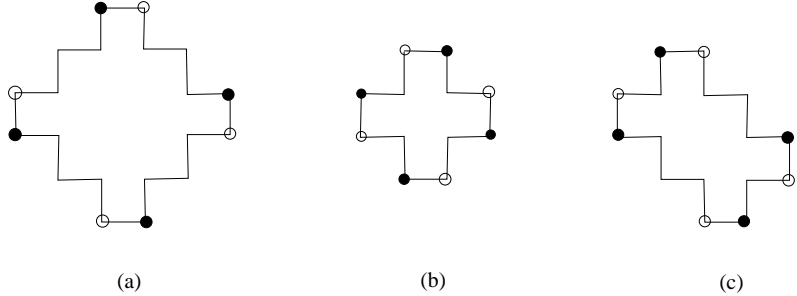}
    \caption{The \emph{diamond tile} (a), the \emph{cross tile} (b) and the \emph{butterfly tile} (c) are three examples of prime double square tiles. Their double BN-factorizations are highlighted on the boundary with dots of different colors.}
    \label{fig:ex}
\end{figure}
\end{Ex}

\begin{Ex}\label{ex:composizione}
    The double square $$P=101\overline{0}\overline{0}101\overline{0}\overline{0}101001010010100\overline{1}\overline{0}\overline{1}00\overline{1}\overline{0}\overline{1}00\overline{1}\overline{0}\overline{1}\overline{0}\overline{0}\overline{1}\overline{0}\overline{1}\overline{0}\overline{0}\overline{1}\overline{0}\overline{1}\overline{0}\overline{0}$$
    is not prime. The trivial morphism, $\varphi_{P}(0)=101\overline{0}\overline{0}101\overline{0}\overline{0}101$, $\varphi_{P}(1)=001010010100$, can be decomposed as $\varphi_{1}(0)=10101$, $\varphi_{1}(1)=0\overline{1}0\overline{1}0$ and $\varphi_{2}(0)=00$, $\varphi_{2}(1)=101$. The composition $\varphi_{P}=\varphi_{1}\circ\varphi_{2}$ maps the simple square in $P$ using the diamond as intermediate step (as shown in  Fig.~\ref{fig:morph}).

    \begin{figure}[ht]
        \centering
        \includegraphics[height=5cm]{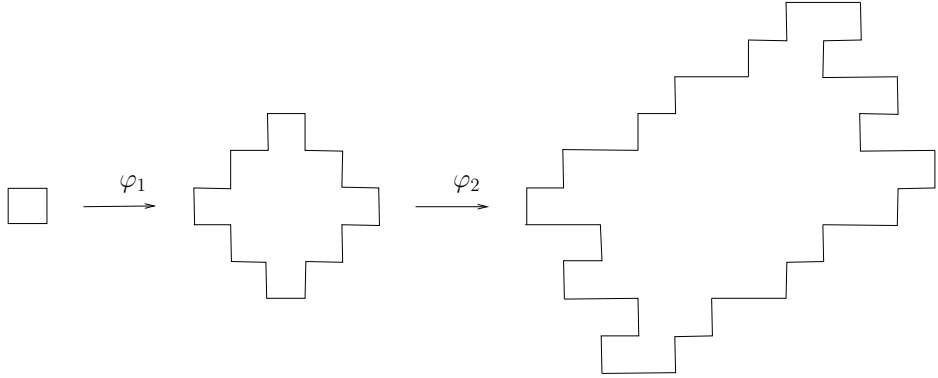}
        \caption{The figure shows the composition of two morphisms that lead to a non-prime double square.}
        \label{fig:morph}
    \end{figure}
\end{Ex}

So, we have introduced all the notions to state a conjecture by Blondin Massé et al. in \cite{BGL}, that constitutes the focus of our work

\medskip
\noindent {\bf Conjecture $\mathbf{35}$.} \emph{Let $P$ be the boundary word of a double square and $\alpha$ a symbol of the alphabet $\Sigma$. If $P$ is prime, then $\alpha\alpha$ is not a subword of $P$.}

\medskip

\emph{Notation:} for each $\alpha \in \Sigma$, we use the notation $\alpha\alpha\subseteq w$ (resp. $\alpha\alpha\not\subseteq w$) to indicate that the word $w$ contains (resp. does not contain) the subword $\alpha\alpha$.

On the other hand, we indicate $w$ to be {\it couple free}, briefly {\it c.f.}, if no two consecutive occurrences of a same letter of $\Sigma$ are present. 

We conclude this section by stating two useful technical lemmas.

\begin{Lemma} [\cite{BGL}]\label{lem:primo}
    Given a double square with BN-factorizations $P=AB\widehat{A}\widehat{B}=XY\widehat{X}\widehat{Y}$, if $P$ is prime then all its factors $A,B,X,Y$ are palindrome.
\end{Lemma}

\begin{Lemma}\label{lemma:palidromi_centrale}
    Given $a,b,u\in\Sigma^{*}$, if $b\Tilde{u}a$ and $b\widehat{u}a$ are palindrome, then $b=\Tilde{a}$ and $u$ is palindrome. If $\Tilde{a}u\Tilde{b}$ and $\overline{b}\Tilde{u}\overline{a}$ are palindrome, then the same result holds.
\end{Lemma}

The proof is trivial.

\section{Properties of the boundary word of prime double squares}\label{sec:properties}

In this section we study the prime double squares' boundary words. In particular, we  exploit some of their properties that, step by step, lead to write them in a  useful general form.

\begin{Lemma}[\cite{BGL}, Lemma~31]\label{lemma:conjugate_wi}
Let $P=w_{1}w_{2}w_{3}w_{4}w_{5}w_{6}w_{7}w_{8}$ be the BN-factorization of (the boundary word of) a prime double square as in (\ref{eq:genform}). It holds that $w_{i+4}=\overline{w}_{i}$ for all $i=1,\dots,4$.
\end{Lemma}




It follows that the boundary word of a prime double square can be written as
\begin{equation}\label{eq:forma2}
P=w_{1}w_{2}w_{3}w_{4}\overline{w}_{1}\overline{w}_{2}\overline{w}_{3}\overline{w}_{4}.
\end{equation}

\begin{Lemma}\label{lemma:alpha_factor}
Let $P=AB\widehat{A}\widehat{B}$ be a prime double square and $\alpha\in\Sigma$. If $\alpha\alpha\subseteq P$, then $\alpha\alpha \subseteq A$ or $\alpha\alpha \subseteq B$. 
\end{Lemma}

\emph{\proofname.}
By contradiction. Let us suppose that the BN-factorization splits $\alpha\alpha$, i.e., there exist $k_{1},k_{2}\in\Sigma^{*}$ such that $P=AB\widehat{A}\widehat{B}=(k_{1}\alpha)(\alpha k_{2})\widehat{A}\widehat{B}$. Since $P$ is a prime double square, its factors are both palindrome by Lemma~\ref{lem:primo}, i.e. $A=(\alpha t_{1}\alpha)$ and $B=(\alpha t_{2}\alpha)$ for some $t_{1},t_{2}\in\Sigma^{*}$. Replacing $A$ and $B$ in $P$, we get

$$P=(\alpha t_{1}\alpha)(\alpha t_{2}\alpha)(\overline{\alpha}\widehat{t}_{1}\overline{\alpha})(\overline{\alpha}\widehat{t}_{2}\overline{\alpha}),$$

reaching a contradiction since the subword $\alpha\overline{\alpha}$ represents a closed path, not allowed in the boundary of a polyomino. \qed

\begin{Cor}\label{cor:split_w}
From the previous lemma and the general form of $P$ in (\ref{eq:forma2}), it directly follows that the term $\alpha\alpha$ can not be split between two consecutive factors $w_i$ and $w_{i+1}$. 
\end{Cor}

\begin{Cor}\label{cor:lunghezze_w}
If $P=w_{1}\dots w_{8}$ is a prime double square, then $|w_{i}|\neq|w_{i+1}|$ for all $i=1,\dots,8$.
\end{Cor}

\emph{\proofname.}
By contradiction. Let be $|w_{i}|=|w_{i+1}|$ for some $i$. Since $w_{i}w_{i+1}$ is a BN-factor and $P$ is prime, then, by Lemma~\ref{lem:primo}, it is a palindrome, i.e., $w_{i}w_{i+1}=w_{i}\Tilde{w_{i}}$. It follows that the last letter of $w_{i}$ and the first one of $w_{i+1}$ match, i.e., we can split a term $\alpha\alpha$ between two BN-factors, in contradiction with Lemma~\ref{lemma:alpha_factor}.\qed

\smallskip

When considering the boundary word $P=w_{1}\dots w_{8}$ of a generic (not necessarily prime) double square, the following property holds:

\begin{Pro}[\cite{BGL}]\label{prop:uv} 
For $i=1,\dots,8$, there exist $u_{i}, v_{i}\in\Sigma^{*}$ and $n_{i}\geq 0$ such that

\begin{equation*}
    \begin{cases}
    
    w_{i}=(u_{i}v_{i})^{n_{i}}u_{i},\\
    \widehat{w}_{i-3}w_{i-1}=u_{i}v_{i}.\\
    \end{cases}
\end{equation*}
\end{Pro}

From Property~\ref{prop:uv} and Lemma~\ref{lemma:conjugate_wi}, we can refine the generic form of the boundary word of a prime double square provided in equation (\ref{eq:forma2}).

First, we note that $w_{i+4}=\overline{w}_{i}$ and $w_{i}=(u_{i}v_{i})^{n_{i}}u_{i}$ 
imply

\begin{equation*}
    \begin{array}{l}
    
    u_{i+4}=\overline{u}_{i},\\
    v_{i+4}=\overline{v}_{i},\\
    n_{i+4}=n_{i},
    
    \end{array}
\end{equation*}

for $i=1,\dots,4$.

By the second equation in Property~\ref{prop:uv}, we can write $u_{i}v_{i}$ depending on $w_{i}$, and replace them in the first one getting 

\begin{equation}\label{eq:parola_bordo}
\begin{split}
P=&(\Tilde{w}_{2}\overline{w}_{4})^{n_{1}}u_{1}(\Tilde{w}_{3}w_{1})^{n_{2}}u_{2}(\Tilde{w}_{4}w_{2})^{n_{3}}u_{3}(\widehat{w}_{1}w_{3})^{n_{4}}u_{4}\dots \\
\dots&(\widehat{w}_{2}w_{4})^{n_{1}}\overline{u}_{1}(\widehat{w}_{3}\overline{w}_{1})^{n_{2}}\overline{u}_{2}(\widehat{w}_{4}\overline{w}_{2})^{n_{3}}\overline{u}_{3}(\Tilde{w}_{1}\overline{w}_{3})^{n_{4}}\overline{u}_{4}.
\end{split}
\end{equation}

The following Lemma is a rephrase of Lemma~11 in \cite{BGL}:

\begin{Lemma}\label{lemma:esponenti}
Given the boundary word $P$ of a prime double square according to equation (\ref{eq:parola_bordo}), there are no two consecutive exponents $n_{1}$, $n_{2}$, $n_{3}$, and $n_{4}$ that are positive.
\end{Lemma}




           
           
           
           
          

\begin{Lemma}
The boundary word $P$ of a prime double square, as in (\ref{eq:parola_bordo}), can always be written in terms of the elements $u_{i}$ only, for $i=1,\dots,4$.
\end{Lemma}

\emph{\proofname.}
If $n_{i}=0$ for all $i$, the thesis trivially holds. So, let us assume w.l.g. $n_{1}>0$. By Lemma~\ref{lemma:esponenti}, it follows $n_{2}=n_{4}=0$, and then $w_{2}=u_{2}$ and $w_{4}=u_{4}$. Since $w_{1}$ and $w_{3}$ are expressed in terms of $w_{2}$ and $w_{4}$ only, the thesis follows. 

\qed 

We finally get the generic form of the boundary word of a prime double square:

\begin{equation}\label{eq:boundary_word}
\begin{split}
    P=&(\Tilde{u}_{2}\overline{u}_{4})^{n_{1}}u_{1}\vdots (\Tilde{u}_{3}u_{1})^{n_{2}}u_{2}\vdots (\Tilde{u}_{4}u_{2})^{n_{3}}u_{3}\vdots (\widehat{u}_{1}u_{3})^{n_{4}}u_{4}|  \\
    &(\widehat{u}_{2}u_{4})^{n_{1}}\overline{u}_{1}\vdots (\widehat{u}_{3}\overline{u}_{1})^{n_{2}}\overline{u}_{2}\vdots (\widehat{u}_{4}\overline{u}_{2})^{n_{3}}\overline{u}_{3}\vdots (\Tilde{u}_{1}\overline{u}_{3})^{n_{4}}\overline{u}_{4}.
    \end{split}
\end{equation}

From now on, we will use a vertical bar to indicate half the boundary word, and dashed, vertical bars to separate two consecutive words $w_{i}$ and $w_{i+1}$.

\begin{Thm}\label{thm:alpha_u}
Let $P$ be the boundary word of a prime double square as in (\ref{eq:boundary_word}), and $\alpha\in\Sigma$. If $\alpha\alpha\subseteq P$, then $\alpha\alpha$ is entirely contained in a factor $u_{i}$, i.e., the term $\alpha\alpha$ can not be split between two consecutive factors $u_{i}$.
\end{Thm}

\emph{\proofname.}
Let us suppose w.l.g. $\alpha\alpha\subseteq w_{1}=(\Tilde{u}_{2}\overline{u}_{4})^{n_{1}}u_{1}$. By contradiction, we suppose to split the occurrences of $\alpha$ between two factors $u_{i}$. We analyze two possible cases for the values of $n_1$.


Let be $n_{1}\geq 2$ and $\alpha\alpha\subseteq\overline{u}_{4}\Tilde{u}_{2}$. Since $\Tilde{u}_{2}$ begins and $\overline{u}_{4}$ ends with $\alpha$, then $\widehat{u}_{2}$ begins with $\overline{\alpha}$ and ${u}_{4}$ ends with $\overline{\alpha}$. In this case, we get

$$w_{4}|w_{5}=\dots u_{4}|\widehat{u}_{2}\dots$$

and $\overline{\alpha}\overline{\alpha}$ is split between two factors $w_{i}$, in contradiction with Lemma~\ref{lemma:alpha_factor}.
    
Let $n_{1}>0$ and $\alpha\alpha\subseteq\Tilde{u}_{2}\overline{u}_{4}$. Since $\Tilde{u}_{2}$ ends and $\overline{u}_{4}$ begins with $\alpha$, then $u_{2}$ begins with $\alpha$ and $u_{4}$ begins with $\overline{\alpha}$. We analyze $w_{3}=(\Tilde{u}_{4}u_{2})^{n_{3}}u_{3}$. If $n_{3}>0$, then  $\overline{\alpha}\alpha\subseteq w_{3}$, and the boundary of the polyomino intersects itself. Then, it holds $n_{3}=0$, and $X=w_{2}w_{3}\stackrel{n_{2}=0}{=}u_{2}u_{3}$. Since $X$ is palindrome, $u_{3}$ ends with $\alpha$, getting $\alpha\overline{\alpha}\subseteq w_{3}w_{4}$, and again the border intersects itself.




A similar reasoning leads to contradiction when considering either  $\alpha\alpha\subseteq\overline{u}_{4}u_{1}$ and $n_{1}>0$ or $\alpha\alpha\subseteq u_{1}\Tilde{u}_{3}$ or $u_{1}u_{2}$ and $n_{1}=0$. 
\qed

\medskip

Since $P$ is a circular word, from now on we suppose, w.l.g., that $|u_{1}|\leq|u_{i}|$ for $i=2,3,4$. Assuming that, it is possible to define $P$ of equation (\ref{eq:boundary_word}) in terms of two factors only, $u_{1}$ and $u_{3}$, and two words $k,p\in\Sigma^{*}$, with $k,p\neq\varepsilon$. It will be clear from the proof of Theorem~\ref{thm:word_u1u3} that $u_2=k\Tilde{u_1}$ and $u_4=\widehat{u}_{1}\overline{p}$. 

Relying on Lemmas~\ref{lem:primo} and~\ref{lemma:esponenti}, we get the following taxonomy according to the positive values of $n_{i}$, with $i=1,\dots,4$:

\begin{description}
    
    \item[a)] $P=u_{1}\vdots k\Tilde{u}_{1}\vdots u_{3} \vdots \widehat{u}_{1}\overline{p}| \overline{u}_{1}\vdots \overline{k}\widehat{u}_{1} \vdots \overline{u}_{3} \vdots \Tilde{u}_{1}p$. This form requires $k$ and $p$ palindrome.

    \item[b)] $P=(u_{1}k\Tilde{u}_{1}p)^{n_{1}}u_{1}\vdots k\Tilde{u}_{1}\vdots u_{3} \vdots \widehat{u}_{1}\overline{p}| (\overline{u}_{1}\overline{k}\widehat{u}_{1}\overline{p})^{n_{1}}\overline{u}_{1}\vdots \overline{k}\widehat{u}_{1} \vdots \overline{u}_{3} \vdots \Tilde{u}_{1}p$. Here $k$ and $p$ are palindrome.

    \item[c)] $P=u_{1}\vdots (\Tilde{u}_{3}u_{1})^{n_{2}}k\Tilde{u}_{1}\vdots u_{3} \vdots \widehat{u}_{1}\overline{p}| \overline{u}_{1}\vdots (\widehat{u}_{3}\overline{u}_{1})^{n_{2}}\overline{k}\widehat{u}_{1} \vdots \overline{u}_{3} \vdots \Tilde{u}_{1}p$. Here $p$ is palindrome.

    \item[d)] $P=u_{1}\vdots k\Tilde{u}_{1}\vdots (\overline{p}\overline{u}_{1}k\Tilde{u}_{1})^{n_{3}}u_{3} \vdots \widehat{u}_{1}\overline{p}| \overline{u}_{1}\vdots \overline{k}\widehat{u}_{1} \vdots (p u_{1}\overline{k}\widehat{u}_{1})^{n_{3}}\overline{u}_{3} \vdots \Tilde{u}_{1}p$. In this case it holds $k$ and $p$ palindrome.

    \item[e)] $P=u_{1}\vdots k\Tilde{u}_{1}\vdots u_{3} \vdots (\widehat{u}_{1}u_{3})^{n_{4}}\widehat{u}_{1}\overline{p}| \overline{u}_{1}\vdots \overline{k}\widehat{u}_{1} \vdots \overline{u}_{3} \vdots (\Tilde{u}_{1}\overline{u}_{3})^{n_{4}}\Tilde{u}_{1}p$. Here $k$ is palindrome.

    \item[f)] $P=(u_{1}k\Tilde{u}_{1}p)^{n_{1}}u_{1}\vdots k\Tilde{u}_{1}\vdots (\overline{p}\overline{u}_{1}k\Tilde{u}_{1})^{n_{3}}u_{3} \vdots \widehat{u}_{1}\overline{p}| (\overline{u}_{1}\overline{k}\widehat{u}_{1}\overline{p})^{n_{1}}\overline{u}_{1}\vdots \overline{k}\widehat{u}_{1} \vdots (p u_{1}\overline{k}\widehat{u}_{1})^{n_{3}}\overline{u}_{3} \vdots \Tilde{u}_{1}p$. Here $k$ and $p$ are palindrome.

    \item[g)] $P=u_{1}\vdots (\Tilde{u}_{3}u_{1})^{n_{2}}k\Tilde{u}_{1}\vdots u_{3} \vdots (\widehat{u}_{1}u_{3})^{n_{4}}\widehat{u}_{1}\overline{p}| \overline{u}_{1}\vdots (\widehat{u}_{3}\overline{u}_{1})^{n_{2}}\overline{k}\widehat{u}_{1} \vdots \overline{u}_{3} \vdots (\Tilde{u}_{1}\overline{u}_{3})^{n_{4}}\Tilde{u}_{1}p$. In this case, $k$ and $p$ are generic, non-empty words in $\Sigma^{*}$.
\end{description}

\begin{Thm}\label{thm:word_u1u3}
   The above seven cases $a) \dots g)$ describe all the possible boundary words of a prime double square.
\end{Thm}

\emph{Sketch of proof.}
$P$ is a prime double square, then its BN-factors $A=w_{1}w_{2}$, $B=w_{3}w_{4}$, $X=w_{2}w_{3}$ and $Y=w_{4}\overline{w}_{1}$ are all palindrome; we further remind that $u_{i}$ is a prefix of $w_{i}$ by Property~\ref{prop:uv}. Then, for any possible value of $n_{i}$ and by the hypothesis on the lengths $|u_{i}|$, it hold $u_{2}=k\Tilde{u}_{1}$ and $u_{4}=\widehat{u}_{1}\overline{p}$ for some $k,p\in\Sigma^{*}$. We notice that we can not deduce an analogous form for the word $u_{3}$, since we do not know the relative lengths w.r.t. $u_{2}$ and $u_{4}$. Analyzing palindromicity properties of the BN-factors, for each case $a) \dots g)$ we deduce the related properties of $k$ and $p$ that lead to the thesis.\qed

\section{Proof of the conjecture}\label{sec:conjecture}

We will prove Conjecture $35$ 
relying on the results of the previous section, in particular the final forms $a) \dots g)$ of the boundary word of a prime double square $P$ (Theorem~\ref{thm:word_u1u3}). We start analyzing case $a)$, where $n_{i}=0$ for all $i$, and then we show how to generalize the results to the remaining cases.


\begin{Lemma}\label{lemma:only_u1}
    Let $P$ be as in case $a)$, and $\alpha \in \Sigma$. If $\alpha\alpha\subseteq u_{1}$, and $u_{3}$, $k$ and $p$ are c.f. (couple free), then it holds that $P$ is not a prime double square.
\end{Lemma}
\emph{\proofname}.
Let us suppose there is one only occurrence $\alpha\alpha\subseteq u_{1}$, i.e., $u_{1}=a\alpha\alpha b$ for some $a,b\in\Sigma^{*}$ and $a$ and $b$ couple free. 
Then
$$P=a\alpha\alpha b\vdots k\Tilde{b}\alpha\alpha\Tilde{a}\vdots u_{3}\vdots \widehat{b}\overline{\alpha}\overline{\alpha}\widehat{a}\overline{p}|\overline{a}\overline{\alpha}\overline{\alpha}\overline{b}\vdots \overline{k}\widehat{b}\overline{\alpha}\overline{\alpha}\widehat{a}\vdots \overline{u}_{3}\vdots \Tilde{b}\alpha\alpha\Tilde{a}p.$$
Since $w_{2}w_{3}$ and $w_{3}w_{4}$ are palindrome, and $\alpha\alpha\not\subseteq a$, $b$, $k$ and $u_{3}$, we get

$$
\begin{cases}
    k\Tilde{b}=\Tilde{u}_{3}a, \\
    u_{3}\widehat{b}=\widehat{p}\overline{a}, \\
    k\Tilde{b}\Tilde{a}u_{3} \text{ palindrome, } \\
    u_{3}\widehat{b}\widehat{a}\overline{p} \text{ palindrome, }
\end{cases}
$$
so simple computations on the factors' lengths lead to $|k|=|p|=|u_{3}|$ and $\Tilde{b}a$ and $u_{3}$ both palindrome. It also follows $|a|=|b|$, with $b=\Tilde{a}$, and $k=\overline{p}=u_{3}$. Replacing in $P$ we get
$$P=a\alpha\alpha\Tilde{a}u_{3}a\alpha\alpha\Tilde{a}u_{3}\overline{a}\overline{\alpha}\overline{\alpha}\widehat{a}u_{3}\overline{a}\overline{\alpha}\overline{\alpha}\widehat{a}\overline{u}_{3}\overline{a}\overline{\alpha}\overline{\alpha}\widehat{a}\overline{u}_{3}a\alpha\alpha\Tilde{a}\overline{u}_{3}.$$

Finally, $P$ can be obtained from the boundary word of the cross, $Q=1010\overline{1}0\overline{1}\overline{0}\overline{1}\overline{0}1\overline{0}$, through the morphism $\varphi(0)=u_{3}$, $\varphi(1)=a\alpha\alpha\Tilde{a}=u_{1}$. It follows that $P$ is not prime.

Let us now conclude the proof by considering the case where $u_{1}$ contains more occurrences of two consecutive equal letters, i.e., $u_{1}=a\alpha\alpha c\beta\beta b$ with $a$ and $b$ c.f. and $c\in\Sigma^{*}$. 
Replacing $u_1$ in $P$ we get
$$P=a\alpha\alpha c\beta\beta b \vdots k\Tilde{b}\beta\beta\Tilde{c}\alpha\alpha\Tilde{a} \vdots u_{3} \vdots \widehat{b}\overline{\beta}\overline{\beta}\widehat{c}\overline{\alpha}\overline{\alpha}\widehat{a} \overline{p}| \dots \Tilde{b}\beta\beta\Tilde{c}\alpha\alpha\Tilde{a}p.
$$

Following the same argument, we immediately argue $\alpha=\beta$ and we get again $u_{3}=\overline{p}=k$ palindrome and $b=\Tilde{a}$; moreover, $c$ is palindrome too. Replacing in $P$, we get the final form 
$$P=a \alpha\alpha c\alpha\alpha\Tilde{a}\vdots u_{3} a\alpha\alpha c\alpha\alpha\Tilde{a}\vdots u_{3}\vdots \overline{a}\overline{\alpha}\overline{\alpha}\overline{c}\overline{\alpha}\overline{\alpha}\widehat{a}u_{3}| \dots  a\alpha\alpha c\alpha\alpha \Tilde{a} \overline{u}_{3},
$$
that can be obtained from the cross through the morphism $\varphi'(0)=u_{3}$, $\varphi'(1)=a\alpha\alpha c\alpha\alpha\Tilde{a}=u_{1}$, that easily generalizes $\varphi$. Even in this case, $P$ is not prime. \qed

\begin{Remark}
    The morphisms $\varphi$ and $\varphi'$ in Lemma~\ref{lemma:only_u1} are well defined and preserve the hat operation, since $u_{1}$ and $u_{3}$ are both palindrome.
\end{Remark}

It is easy to verify that, as shown in the proof of Lemma~\ref{lemma:only_u1}, when considering more occurrences of couples of consecutive equal letters in a subword of $P$, the morphisms, if any, that map a prime double square into $P$ are a simple generalization of the morphism obtained by considering one occurrence of two equal letters only. For this reason, the next results will be provided adding this further hypothesis. 

\begin{Lemma}\label{lemma:p_not_u1}
Let $P$ be a double square as in case $a)$, and $\alpha \in \Sigma$. $P$ is not prime when one of the following conditions verifies:
\begin{itemize}
    \item[i)] $u_{1}$ is c.f. and $\alpha\alpha\subseteq p$; 
    \item[ii)] $u_{1}$ is c.f. and $\alpha\alpha\subseteq k$;
    \item[iii)] $u_{1}$ is c.f. and $\alpha\alpha\subseteq u_{3}$. 
\end{itemize}
\end{Lemma}

\emph{\proofname}.
\noindent $i)$ Since $P$ has the form $a)$, then each exponent $n_{i}$ in equation (\ref{eq:boundary_word}) equals $0$, and $p$ can be written as the palindrome $a\alpha\alpha\Tilde{a}$. 
From $u_{1},a$ couple free, and $w_{3}w_{4}$ palindrome, it follows  $u_{3}=\overline{a}\overline{\alpha}\overline{\alpha}\Tilde{y}$ for suitable $x,y\in\Sigma^{*}$ s.t. $\overline{a}=xy$. Then
$$P=u_{1}\vdots k\Tilde{u}_{1}\vdots xy\overline{\alpha}\overline{\alpha}\Tilde{y} \vdots \widehat{u}_{1}xy\overline{\alpha}\overline{\alpha}\Tilde{y}\Tilde{x}|\overline{u}_{1}\vdots \overline{k}\widehat{u}_{1}\vdots \overline{x}\overline{y}\alpha\alpha\widehat{y} \vdots \Tilde{u}_{1}\overline{x}\overline{y}\alpha\alpha\widehat{y}\widehat{x}.$$

By the palindromicity of $w_{2}w_{3}$, we argue that $k=y\overline{\alpha}\overline{\alpha}t$ for some $t\in\Sigma^{*}$. We conclude that $t=\Tilde{y}$, since $k$ is palindrome by hypothesis ($n_{i}=0$ for all $i$). Replacing again in $P$, we get

$$P=u_{1}\vdots y\overline{\alpha}\overline{\alpha}\Tilde{y}\Tilde{u}_{1}\vdots xy\overline{\alpha}\overline{\alpha}\Tilde{y}\vdots \widehat{u}_{1}xy\overline{\alpha}\overline{\alpha}\Tilde{y}\Tilde{x}|\overline{u}_{1}\vdots \overline{y}\alpha\alpha\widehat{y}\widehat{u}_{1}\vdots \overline{x}\overline{y}\alpha\alpha\widehat{y}\vdots \Tilde{u}_{1}\overline{x}\overline{y}\alpha\alpha\widehat{y}\widehat{x}$$

and, by the palindromicity of $w_{2}w_{3}$ and $w_{3}w_{4}$, we observe that $\Tilde{u}_{1}x$ and $\widehat{u}_{1}x$ are both palindrome, and then $x=\varepsilon$ and $u_{1}=\Tilde{u}_{1}$. We get the final form

$$P=u_{1}\vdots y\overline{\alpha}\overline{\alpha}\Tilde{y}{u}_{1}\vdots y \overline{\alpha}\overline{\alpha}\Tilde{y}\vdots \overline{u}_{1}y\overline{\alpha}\overline{\alpha}\Tilde{y}|\overline{u}_{1}\vdots \overline{y}\alpha\alpha\widehat{y}\overline{u}_{1}\vdots \overline{y}\alpha\alpha\widehat{y}\vdots u_{1}\overline{y}\alpha\alpha\widehat{y}.$$

A morphism can now be defined mapping the cross to $P$, i.e., $\varphi(0)=y\overline{\alpha}\overline{\alpha}\Tilde{y}=\overline{a}\overline{\alpha}\overline{\alpha}\widehat{a}=u_{3}$, $\varphi(1)=u_{1}$. Then $P$ is not prime.

\noindent The proofs of cases $ii)$ and $iii)$ are similar to $i)$.\qed

\begin{Cor}
    Let $P$ be a double square as in case $a)$, and $\alpha \in \Sigma$. If $\alpha\alpha\subseteq P$ and $u_{1}$ is c.f., then $\alpha\alpha\subseteq p,k,u_{3}$.
\end{Cor}

\begin{Lemma}\label{lemma:only_u1_p}
Let $P$ be a double square as in case $a)$. $P$ is not prime when one of the following conditions verifies:
\begin{itemize}
    \item[i)] $u_{1}$ and $p$ are not c.f. while  $u_{3}$ and $k$ are;
    \item[ii)] $u_{1}$ and $k$ are not c.f. while  $u_{3}$ and $p$ are; 
    \item[iii)] $u_{1}$, $p$ and $k$ are not c.f. while  $u_{3}$ is.
\end{itemize}
\end{Lemma}
\emph{\proofname}.
$i)$ Let us assume that one only occurrence of a couple of consecutive equal letters is both in $u_1$ and $p$, say $u_{1}=a\alpha\alpha b$ and $p=c\beta\beta\Tilde{c}$ (recall that $p$ is palindrome), with $\alpha,\beta \in \Sigma$, and $a$, $b$, $c$ that are couple free. 


Since $w_{2}w_{3}$ and $w_{3}w_{4}$ are palindrome and $\alpha\alpha\not\subseteq a,b,c,k,u_{3}$, we deduce $\alpha=\beta$, $k\Tilde{b}=\Tilde{u}_{3}a$ and $\widehat{c}=\overline{b}\Tilde{u}_{3}$. Replacing in $P$ we obtain
$$P=a\alpha\alpha b\vdots\Tilde{u}_{3}a\alpha\alpha\Tilde{a}\vdots u_{3}\vdots\widehat{b}\overline{\alpha}\overline{\alpha}\widehat{a}u_{3}\widehat{b}\overline{\alpha}\overline{\alpha}\overline{b}\Tilde{u}_{3}|\overline{a}\overline{\alpha}\overline{\alpha}\overline{b}\vdots \widehat{u}_{3}\overline{a}\overline{\alpha}\overline{\alpha}\widehat{a}\vdots\overline{u}_{3}\vdots\Tilde{b}\alpha\alpha\Tilde{a}\overline{u}_{3}\Tilde{b}\alpha\alpha b\widehat{u}_{3}.$$

By the palindromicity of $w_{1}w_{2}$ and $w_{3}w_{4}$, it follows that also $b\Tilde{u}_{3}{a}$ and $\widehat{a}u_{3}\widehat{b}$ are. Since $b=\Tilde{a}$ and $u_{3}=\Tilde{u}_{3}$ by Lemma~\ref{lemma:palidromi_centrale}, then $P$ can be expressed as 
$$P=a\alpha\alpha\Tilde{a}\vdots u_{3}a\alpha\alpha\Tilde{a}\vdots u_{3}\vdots\overline{a}\overline{\alpha}\overline{\alpha}\widehat{a}u_{3}\overline{a}\overline{\alpha}\overline{\alpha}\widehat{a}u_{3}|\overline{a}\overline{\alpha}\overline{\alpha}\widehat{a}\vdots\overline{u}_{3}\overline{a}\overline{\alpha}\overline{\alpha}\widehat{a}\vdots\overline{u}_{3}\vdots a\alpha\alpha\Tilde{a}\overline{u}_{3}a\alpha\alpha\Tilde{a}\overline{u}_{3}.$$
The morphism $\varphi(0)=u_{3}$, $\varphi(1)=a\alpha\alpha\Tilde{a}=u_{1}$ maps the butterfly $Q=1010\overline{1}0\overline{1}0\overline{1}\overline{0}\overline{1}\overline{0}1\overline{0}1\overline{0}$ into $P$, preventing it from being prime.

As already observed, the cases of more occurrences of couples of the same letters in $u_1$ and $p$ can be treated similarly.

\noindent Also the proofs of cases $ii)$ and $iii)$ are similar.\qed

\begin{Remark}\label{rem:couple}
    Since the morphisms we defined in the proofs of Lemmas~
    \ref{lemma:p_not_u1} and~\ref{lemma:only_u1_p} map $0,1$ in $u_{1},u_{3}$, and $k,p$ result to be equal, unless conjugation, to $u_{1}$ or $u_{3}$ or their concatenation (as in Lemma~\ref{lemma:only_u1_p}), we can generalize these results when more couples of consecutive equal letters occur in one or more subwords of $P$, as shown in the proof of Lemma~\ref{lemma:only_u1}.
\end{Remark}

The analysis of the possible positions in $P$ of the subword $\alpha \alpha $ continues, and two more lemmas are provided. We choose to set them separately since the proofs, even though both proceeding by contradiction, are different. 

\begin{Lemma}\label{lemma:only_u1_u3}
 Let $P$ be a double square as in case $a)$. If $u_{1}$ and $u_{3}$ are not c.f., while $k$ and $p$ are, then $P$ is not prime.
\end{Lemma}

\emph{\proofname}.
 Let us assume that one only occurrence of a couple of consecutive equal letters is both in $u_1$ and $u_3$, say $u_{1}=a\alpha\alpha b$ and $u_{3}=c\beta\beta d$, with $\alpha,\beta \in \Sigma$, and $a$, $b$, $c$, $d$ that are couple free.
Replacing in $P$, we get

$$P=a\alpha\alpha b\vdots k\Tilde{b}\alpha\alpha\Tilde{a} \vdots c\beta\beta d \vdots \widehat{b}\overline{\alpha}\overline{\alpha}\widehat{a}\overline{p} | \overline{a}\overline{\alpha}\overline{\alpha}\overline{b} \vdots \overline{k}\widehat{b}\overline{\alpha}\overline{\alpha}\widehat{a} \vdots \overline{c}\overline{\beta}\overline{\beta}\overline{d} \vdots \Tilde{b}\alpha\alpha\Tilde{a}p.$$

Looking at the palindrome $w_{2}w_{3}$ and $w_{3}w_{4}$, we deduce (respectively) $\beta=\alpha$ and $\beta=\overline{\alpha}$, that give a contradiction. Then, a double square $P$ containing such configurations of couples $\alpha\alpha$ is not prime. If more than one occurrence of couples of letters are in $u_1$ or $u_3$ a similar proof holds.\qed

\begin{Lemma}\label{lemma:u1_u3_k}
The following statements hold:
\begin{itemize}
    \item[i)] if $u_{1}$, $u_{3}$ and $k$ are not c.f., while  $p$ is, then $P$ is not a prime double square;
    \item[ii)] similarly, if $u_{1}$, $u_{3}$ and $p$ are not c.f., while  $k$ is, then $P$ is not a prime double square.
\end{itemize}
    
\end{Lemma}

\emph{\proofname}.
$i)$ Let us again assume that one only occurrence of a couple of consecutive equal letters is both in $u_1$, $u_3$ and $k$, say  $u_{1}=a\alpha\alpha b$, $u_{3}=c\beta\beta d$ and $k=e\lambda\lambda\Tilde{e}$,  with $\alpha,\beta,\lambda \in \Sigma$, and $a$, $b$, $c$, $d$, $e$ that are couple free.
From the palindromicity of $w_{2}w_{3}$ and $w_{3}w_{4}$, we have $\lambda=\beta=\overline{\alpha}$, $e=\Tilde{d}$, $c=\overline{p}\overline{a}$, $\Tilde{e}\Tilde{b}=\Tilde{c}a$. Moreover, we highlight that $d\widehat{b}$ is palindrome. So $P$ can be written as 
$$P=a\alpha\alpha b\vdots \Tilde{d}\overline{\alpha}\overline{\alpha}d\Tilde{b}\alpha\alpha\Tilde{a} \vdots \overline{p}\overline{a}\overline{\alpha}\overline{\alpha}d \vdots \widehat{b}\overline{\alpha}\overline{\alpha}\widehat{a}\overline{p}| \overline{a}\overline{\alpha}\overline{\alpha}\overline{b} \vdots \widehat{d}\alpha\alpha\overline{d}\widehat{b}\overline{\alpha}\overline{\alpha}\widehat{a} \vdots pa\alpha\alpha \overline{d} \vdots \Tilde{b}\alpha\alpha\Tilde{a}p.$$

Similarly, $d\Tilde{b}=\widehat{a}\overline{p}a$ from $w_{2}w_{3}$, and
$$P=a\alpha\alpha\Tilde{a}\overline{p}\overline{a}\overline{\alpha}\overline{\alpha}\widehat{a}\overline{p}a\alpha\alpha\Tilde{a}\overline{p}\overline{a}\overline{\alpha}\overline{\alpha}d\widehat{b}\overline{\alpha}\overline{\alpha}\widehat{a}\overline{p}| \overline{a}\overline{\alpha}\overline{\alpha}\widehat{a}pa\alpha\alpha\Tilde{a}p\overline{a}\overline{\alpha}\overline{\alpha}\widehat{a}pa\alpha\alpha\overline{d}\Tilde{b}\alpha\alpha\Tilde{a}p,$$
with $d\widehat{b}$ palindrome and $d\Tilde{b}=\widehat{a}\overline{p}a$. 

We underline that the subword $w=\overline{p}\overline{a}\overline{\alpha}\overline{\alpha}\widehat{a}pa\alpha\alpha\Tilde{a}$ is such that $|w|_{1}=|w|_{\overline{1}}$ and $|w|_{0}=|w|_{\overline{0}}$; we then conclude that the boundary of the polyomino intersects itself, reaching a contradiction.


The proof of $ii)$ is  similar to $i)$. 
\qed

\begin{Lemma}\label{lemma:all_words}
    If $u_{1}$, $u_{3}$, $p$ and $k$ are not c.f., then $P$ is not a prime double square.
\end{Lemma}

The result can be obtained providing a morphism $\varphi$ that maps the cross in $P$, following the same argument used in the previous proofs. 
We underline again that similar proofs lead to analogous morphisms when considering more couples of consecutive equal letters in $u_{1},u_{3},k$ or $p$.

\begin{Thm}\label{thm:proof_case1}
    If $P$ is the boundary word of a prime double square s.t. $n_{i}=0$ for all $i$, then $P$ is couple free.
\end{Thm}

The proof directly follows from the previous lemmas, where a complete analysis of the positions of consecutive equal letters inside the form $a)$ of $P$ is carried on.

\subsection*{How to generalize the previous results to the remaining forms $b) \dots g)$ of the word $P$}

We indicate here, for each form of the boundary word of a prime double square from $a)$ to $g)$, the BN-factors involved in the lemmas leading to the proof of Theorem~\ref{thm:proof_case1}. As a matter of fact, it is possible to follow the lemmas' sequence   according to the new indicated factors to obtain the related version of Theorem~\ref{thm:proof_case1} for each remaining case from $b)$ to $g)$. Hence, Conjecture 35 is proved. The factors are the following:

\begin{description}
    \item[Case a)] $A=u_{1}k\Tilde{u}_{1}$, $B=u_{3}\widehat{u}_{1}\overline{p}$, $X=k\Tilde{u}_{1}u_{3}$, $Y=\widehat{u}_{1}\overline{p}\overline{u}_{1}$.
    
    \item[Case b)] $A=(u_{1}k\Tilde{u}_{1}p)^{n_{1}}u_{1}k\Tilde{u}_{1}$, $B=u_{3}\widehat{u}_{1}\overline{p}$, $X=k\Tilde{u}_{1}u_{3}$, $Y=\widehat{u}_{1}\overline{p}(\overline{u}_{1}\overline{k}\widehat{u}_{1}\overline{p})^{n_{1}}\overline{u}_{1}$.
   
   \item[Case c)] $A=u_{1}(\Tilde{u}_{3}u_{1})^{n_{2}}k\Tilde{u}_{1}$, $B=u_{3}\widehat{u}_{1}\overline{p}$, $X=(\Tilde{u}_{3}u_{1})^{n_{2}}k\Tilde{u}_{1}u_{3}$, $Y=\widehat{u}_{1}\overline{p}\overline{u}_{1}$.
    
    \item[Case d)] $A=u_{1}k\Tilde{u}_{1}$, $B=(\overline{p}\overline{u}_{1}k\Tilde{u}_{1})^{n_{3}}u_{3}\widehat{u}_{1}\overline{p}$, $X=k\Tilde{u}_{1}(\overline{p}\overline{u}_{1}k\Tilde{u}_{1})^{n_{3}}u_{3}$, $Y=\widehat{u}_{1}\overline{p}\overline{u}_{1}$.
    
    \item[Case e)] $A=u_{1}k\Tilde{u}_{1}$, $B=u_{3}(\widehat{u}_{1}u_{3})^{n_{4}}\widehat{u}_{1}\overline{p}$, $X=k\Tilde{u}_{1}u_{3}$, $Y=(\widehat{u}_{1}u_{3})^{n_{4}}\widehat{u}_{1}\overline{p}\overline{u}_{1}$.
    
    \item[Case f)] $A=(u_{1}k\Tilde{u}_{1}p)^{n_{1}}u_{1}k\Tilde{u}_{1}$, $B=(\overline{p}\overline{u}_{1}k\Tilde{u}_{1})^{n_{3}}u_{3}\widehat{u}_{1}\overline{p}$, $X=k\Tilde{u}_{1}(\overline{p}\overline{u}_{1}k\Tilde{u}_{1})^{n_{3}}u_{3}$, $Y=\widehat{u}_{1}\overline{p}(\overline{u}_{1}\overline{k}\widehat{u}_{1}\overline{p})^{n_{1}}\overline{u}_{1}$.
    
    \item[Case g)] $A=u_{1}(\Tilde{u}_{3}u_{1})^{n_{2}}k\Tilde{u}_{1}$, $B=u_{3}(\widehat{u}_{1}u_{3})^{n_{4}}\widehat{u}_{1}\overline{p}$, $X=(\Tilde{u}_{3}u_{1})^{n_{2}}k\Tilde{u}_{1}u_{3}$, $Y=(\widehat{u}_{1}u_{3})^{n_{4}}\widehat{u}_{1}\overline{p}\overline{u}_{1}$.
\end{description}

Focusing on case a), we notice that $A$ and $Y$ are palindrome by construction, while $B$ and $X$ begin or end with the words $k,p$ or $u_{3}$. As seen in all proofs of the previous lemmas, when including a couple $\alpha\alpha$ in $P$ we compare $u_{3},p$ and $k$ by the palindromicity of $B$ and $X$. Through this procedure, we can write $k$ and $p$ as a concatenation of $u_{3}$ and $u_{1}$ (as highlighted in Remark~\ref{rem:couple}). Then, through the palindromicity of the BN-factors $A$ and $Y$, we deduce $\Tilde{u}_{1}=u_{1}$ and $\Tilde{u}_{3}=u_{3}$, and finally provide the morphism $\varphi$. 

We can generalize this procedure analyzing the BN-factors in the other cases:

\begin{itemize}
    \item [b)] $B$ and $X$ allow to compare $k,p$ and $u_{3}$, while $A$ and $Y$ give the palindromicity of $u_{1},u_{3}$.

    \item [c)] $A$ and $B$ allow to compare $k,p$ and $u_{3}$, while $X$ and $Y$ give the palindromicity of $u_{1},u_{3}$.

    \item [d)] $B$ and $X$ allow to compare $k,p$ and $u_{3}$, while $A$ and $Y$ give the palindromicity of $u_{1},u_{3}$.

    \item [e)] $B,X$ and $Y$ allow to compare $k,p$ and $u_{3}$, while $A$ gives the palindromicity of $u_{1},u_{3}$.

    \item [f)] $B$ and $X$ allow to compare $k,p$ and $u_{3}$, while $A$ and $Y$ give the palindromicity of $u_{1},u_{3}$.

    \item [g)] $A,B$ and $Y$ allow to compare $k,p$ and $u_{3}$, while $X$ gives the palindromicity of $u_{1},u_{3}$.
\end{itemize}

After these observations, we can conclude that the results obtained in case a) admit a generalization to all the remaining cases. So, we can state the final theorem and solve Conjecture 35:

\begin{Thm}\label{thm:proof_conjecture}
If $P$ is the boundary word of a prime double square and $\alpha\in\Sigma$, then $\alpha\alpha\not\subseteq P$.
\end{Thm}

On the other hand, the following example is related to case $g)$ when $u_1$ is not c.f., while $u_{3},p$ and $k$ are. In this case we are able, similarly to case $a)$, to define a morphism that leads to the non primality of $P$.
\begin{Ex}
Let us consider $P$ having the form $g)$. We assume that $u_1$ is not c.f., while $u_{3},p$ and $k$ are. We can map in $P$ the double square $Q=1(01)^{n_{2}}010(\overline{1}0)^{n_{4}}\overline{1}0\overline{1}(\overline{0}\overline{1})^{n_{2}}\overline{0}\overline{1}\overline{0}(1\overline{0})^{n_{4}}1\overline{0}$, depicted in Fig~\ref{fig:ex_gen}, through the morphism $\varphi(0)=u_{3}=k=\overline{p}$, $\varphi_{1}=u_{1}$. The proof is similar to Lemma~\ref{lemma:only_u1}. 

\begin{figure}[ht]
    \centering
    \includegraphics[height=5.5cm]{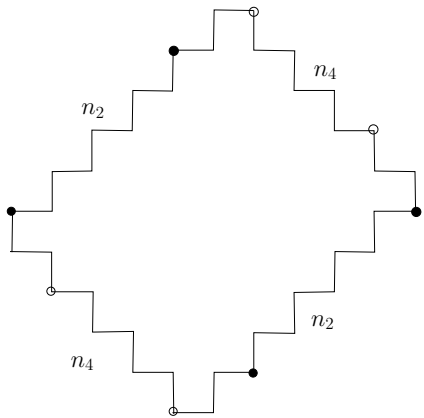}
    \caption{The polyomino $Q$ can be used as an intermediate step to define the trivial morphism from the simple square to a double square $P$ in form g), with $n_{2}=4$, $n_{4}=3$ and s.t. $\alpha\alpha\subseteq u_{1}$ only. Notice that varying the values $n_{2}$ and $n_{4}$ is equivalent to extend or reduce the length of the sides of the polyomino (highlighted with dots), in this case w.r.t. the butterfly.}
    \label{fig:ex_gen}
\end{figure}

\end{Ex}

\section{Conclusion}\label{sec:conclusion}

In this paper we consider the Conjecture 35 of \cite{BGL} and we prove it by showing that there are no couples of the same symbol in the boundary word of a prime double square polyomino. The study of double squares  greatly benefits from the characterization of their prime elements. In particular, it constitutes a step toward the definition of an algorithm that directly generates them, and successively all the double squares, without repetitions or pruning steps. The final aim is to enhance the performances of the algorithm defined in \cite{BGL} by reducing the size of the explored space. As a matter of fact, it could also be interesting moving the spot to the pseudo-hexagon polyominoes and performing a similar inspection. 

From a combinatorial perspective, Theorem~\ref{thm:proof_conjecture} may lead to the definition of a general form of the boundary of a prime double square as a word on a four letters alphabet. Usually, similar expressions lead to the definition of a  growth law according to some parameters, usually perimeter or area, and consequently to the possibility of their enumeration.

\end{document}